\documentclass[12pt]{article}
\newtheorem{thm}{Theorem}[section]
\newtheorem{prop}[thm]{Proposition}
\newtheorem{lem}[thm]{Lemma}

\newtheorem{con}[thm]{Conjecture}

\newcommand{\qed}{\hfill \fbox{}\medskip}
\newcommand{\proof}{\medskip\noindent{\bf Proof.}\quad }

\title{A Conjecture about Raising Operators for Macdonald Polynomials}

\author{Jun'ichi Shiraishi\\
\\
{\it Graduate School of 
Mathematical Science, }\\ {\it University of Tokyo, Tokyo, Japan}}

\date{}

\begin{document}

\maketitle

\maketitle

\begin{abstract}
A multivariable hypergeometric-type formula for raising operators
of the Macdonald polynomials is conjectured. 
It is proved that this agrees with 
Jing and J\'ozefiak's expression for the two-row Macdonald polynomials,
and also with 
Lassalle and Schlosser's formula for partitions with length  three.
\end{abstract}

\section{Introduction}

In this article, we present an observation that
the raising operators for the Macdonald polynomials 
$Q_\lambda(x;q,t)$ \cite{Mac}
can be written in a form of 
multivariable basic hypergeometric-type series.

In the work by Lassalle and Schlosser \cite{LS} (see also \cite{LS0}\cite{L2}), 
the fully explicit formula for the raising operator for the Macdonald 
polynomials was obtained (Theorem 5.1 of \cite{LS}).
It was derived by inverting the Pieri formula for 
the Macdonald polynomial.
Jing and J\'ozefiak's expression 
for the two-row Macdonald polynomials \cite{JJ} is recovered from their 
general formula.
The case when the indexing partition is length three,
was studied in the preceding work by Lassalle \cite{L}.

In the papers \cite{S1,S2}, it was observed that 
a certain class of $n$-fold integral transformations $\{I(\alpha)|\alpha\in {\bf C}\}$ 
forms a commutative family, namely $[I(\alpha),I(\beta)]=0$.
The commutativity was proved only for the simplest case $n=2$
by using the explicit formulas for the eigenfunctions of $I(\alpha)$,
and several summation and transformation formulas for 
the basic hypergeometric series \cite{S1}.
To prove the commutativity for $n\geq 3$ remains as
an open problem, since properties of the eigenfunctions 
have not been studied well. 
It was observed that a modified Macdonald difference operator 
$D^1(s_1,\cdots,s_n,q,t)$ (see (\ref{dif}) below) and the integral transformation $I(\alpha)$
are also commutative with each other. The commutativity was proved for the 
simplest case $n=2$. In Appendix of \cite{S2}, it was shown that 
the eigenfunction of $D^1(s_1,\cdots,s_n,q,t)$ can be interpreted as
a raising operator for the Macdonald polynomials.

An explicit formula for the eigenfunctions of $I(\alpha)$ or 
$D^1(s_1,\cdots,s_n,q,t)$ was conjectured for $n=3$ in \cite{S1}.
(See (\ref{g-fun}) below.) The
 structure of the series (\ref{g-fun}) looks very much different from 
the one obtained by Lassalle and Schlosser.
Therefore, an explanation, which connects
these, is in order. At present, this relation is still unclear.
One may, however, observe that 
Lassalle and Schlosser's formula can be recast in a form, 
which is somewhat closer to the series (\ref{g-fun}).
The aim of this paper is 
to present our observation about this.
\medskip

The plan of the paper is as follows. In Section 2, a conjecture
for the explicit formula of the eigenfunction of $D^1(s_1,\cdots,s_n,q,t)$
is given. The conjecture is recast in the form of the raising operator for
the Macdonald polynomials in Section 3.
Section 4 is devoted to recalling Lassalle and Schlosser's theorem for
the Macdonald polynomials. Then our conjecture is compared with 
Lassalle and Schlosser's result. In Section 5, the case $n=2$ is treated
and our conjecture is proved. In Section 6, it is proved that 
our conjecture for the raising operator agrees with 
Lassalle and Schlosser's formula for $n=3$.
Some special cases $t=q,q^2,q^3,\cdots$, and $q=0$ are discussed in Section 7.

\section{Basic Hypergeometric-like Series}
Let $n$ be a positive integer, and 
$s_1,s_2,\cdots,s_n$ be indeterminates.
Introduce $c_n(\{\theta_{i,j};1\leq i<j\leq n\};s_1,\cdots,s_n,q,t)$
recursively by $c_1(-;s_1,q,t)=1$ and
\begin{eqnarray}
&&c_n(\{\theta_{i,j};1\leq i<j\leq n\};s_1,\cdots,s_n,q,t) \nonumber\\
&=&
c_{n-1}(\{\theta_{i,j};1\leq i<j\leq n-1\};
q^{\theta_{1,n}}s_1,\cdots,q^{\theta_{n-1,n}}s_{n-1},q,t)\label{c_n}\\
&&\times
\prod_{k=1}^{n-1}
t^{\theta_{k,n}} {(qt^{-1};q)_{\theta_{k,n}} \over (q;q)_{\theta_{k,n}} }
{(qt^{-1}s_k/s_n;q)_{\theta_{k,n}} \over (qs_k/s_n;q)_{\theta_{k,n}} }\nonumber\\
&&\times
\prod_{1\leq \ell<m\leq n-1}
 {(qt^{-1} s_\ell/s_m;q)_{\theta_{\ell,n}} \over (q s_\ell/s_m;q)_{\theta_{\ell,n}} }
  {(q^{-\theta_{m,n}}t s_\ell/s_m;q)_{\theta_{\ell,n}} \over 
  (q^{-\theta_{m,n}}s_\ell/s_m;q)_{\theta_{\ell,n}} },
  \nonumber
\end{eqnarray}
where we have used the $q$-shifted factorial 
$(a;q)_n=(1-a )(1-a q)\cdots (1-a q^{n-1})$.
For example, we have
\begin{eqnarray}
&&c_2(\theta_{1,2};s_1,s_2,q,t)=
t^{\theta_{1,2}} {(qt^{-1};q)_{\theta_{1,2}} \over (q;q)_{\theta_{1,2}} }
{(qt^{-1}s_1/s_2;q)_{\theta_{1,2}} \over (qs_1/s_2;q)_{\theta_{1,2}} },\label{c2}\\
&&c_3(\theta_{1,2},\theta_{1,3}, \theta_{2,3};s_1,s_2,s_3,q,t)\nonumber\\
&&=
t^{\theta_{1,2}} {(qt^{-1};q)_{\theta_{1,2}} \over (q;q)_{\theta_{1,2}} }
{(q^{\theta_{1,3}-\theta_{2,3}}qt^{-1}s_1/s_2;q)_{\theta_{1,2}} \over 
(q^{\theta_{1,3}-\theta_{2,3}}qs_1/s_2;q)_{\theta_{1,2}} }\label{c3}\\
&&\times
t^{\theta_{1,3}} {(qt^{-1};q)_{\theta_{1,3}} \over (q;q)_{\theta_{1,3}} }
{(qt^{-1}s_1/s_3;q)_{\theta_{1,3}} \over 
(qs_1/s_3;q)_{\theta_{1,3}} }
t^{\theta_{2,3}} {(qt^{-1};q)_{\theta_{2,3}} \over (q;q)_{\theta_{2,3}} }
{(qt^{-1}s_2/s_3;q)_{\theta_{2,3}} \over 
(qs_2/s_3;q)_{\theta_{2,3}} }\nonumber\\
&&\times
{(qt^{-1}s_1/s_2;q)_{\theta_{1,3}} \over 
(qs_1/s_2;q)_{\theta_{1,3}} }
{(q^{-\theta_{2,3}}ts_1/s_2;q)_{\theta_{1,3}} \over 
(q^{-\theta_{2,3}}s_1/s_2;q)_{\theta_{1,3}} },\nonumber
\end{eqnarray}
and so on. 
The product expression for $c_n(\{\theta_{i,j}\}_{1\leq i<j\leq n};s_1,\cdots,s_n)$
reads as follows
\begin{eqnarray}
&&c_n(\{\theta_{i,j};1\leq i<j\leq n\};s_1,\cdots,s_n,q,t) \nonumber\\
&=&
\prod_{1\leq i<j\leq n}
t^{\theta_{i,j}}
{(qt^{-1};q)_{\theta_{i,j}}\over (q;q)_{\theta_{i,j}} }
{(q^{\sum_{a=j+1}^n (\theta_{i,a}-\theta_{j,a})}qt^{-1}s_i/s_j;q)_{\theta_{i,j}}\over 
(q^{\sum_{a=j+1}^n (\theta_{i,a}-\theta_{j,a})}q s_i/s_j;q)_{\theta_{i,j}} }\\
&&\times
\prod_{k=3}^n
\prod_{1\leq \ell <m\leq k-1}
{(q^{\sum_{b=k+1}^n (\theta_{\ell,b}-\theta_{m,b})}qt^{-1}s_\ell/s_m;q)_{\theta_{\ell,k}}\over 
(q^{\sum_{b=k+1}^n (\theta_{\ell,b}-\theta_{m,b})}q s_\ell/s_m;q)_{\theta_{\ell,k}} }
\nonumber\\
&&\qquad \times
{(q^{-\theta_{m,k}}
q^{\sum_{b=k+1}^n (\theta_{\ell,b}-\theta_{m,b})}t s_\ell/s_m;q)_{\theta_{\ell,k}}\over 
(q^{-\theta_{m,k}}
q^{\sum_{b=k+1}^n (\theta_{\ell,b}-\theta_{m,b})} s_\ell/s_m;q)_{\theta_{\ell,k}} }.
\nonumber
\end{eqnarray}
\medskip

In the paper \cite{S2}, a modified Macdonald difference operator
acting on the space of formal power series 
$F[[x_2/x_1,x_3/x_2,\cdots,x_n/x_{n-1}]]$
was investigated, where $F={\bf Q}(q,t,s_1,s_2,\cdots,s_n)$.
It is defined by
\begin{eqnarray}
&&D^1(s_1,s_2,\cdots,s_n,q,t)=
\sum_{i=1}^n s_i
\prod_{j<i}
{1-q^{- 1} t x_i/x_j\over 1-q^{-1}x_i/x_j}
\prod_{k>i}
{1-q t^{- 1} x_k/x_i\over 1-qx_k/x_i}
T_{q^{-1},x_i},\label{dif}
\end{eqnarray}
where $T_{q,x_i}$ denotes the $q$-shift operator
$T_{q,x_i}  \cdot g(x_1,\cdots,x_n)=
g(x_1,\cdots,q x_i,\cdots,x_n)$
and the rational factors in (\ref{dif}) should be understood as the series
\begin{eqnarray}
{1-q^{\pm 1} t^{\mp 1}x\over 1-q^{\pm 1}x}=
1+\sum_{n=1}^\infty (1-t^{\mp 1}) q^{\pm n}x^n. \label{thetapm}
\end{eqnarray}

Let us consider a basic hypergeometric-like series
\begin{eqnarray}
&&f(x_1,x_2,\cdots,x_n;s_1,\cdots,s_n,q,t)
=\prod_{1\leq k<\ell\leq n}(1-x_\ell/x_k)\label{f}\\
&&\qquad\times
\sum_{(\theta)\in \mathsf{M}^{(n)}}
c_n(\{\theta_{i,j};1\leq i<j\leq n\};s_1,\cdots,s_n,q,t)
\prod_{1\leq i<j\leq n}
(x_j/x_i)^{\theta_{i,j}}, \nonumber
\end{eqnarray}
where $(\theta)=(\theta_{i,j})\in \mathsf{M}^{(n)}$. 
Here we have used the notation introduced in \cite{LS}, 
namely $\mathsf{M}^{(n)}$ denotes the set
of upper triangular $n \times n$
matrices with nonnegative integers, and $0$ on the 
diagonal. 
\medskip

Then we have the following observation.
\begin{con}\label{con}
The series $f(x_1,x_2,\cdots,x_n;s_1,\cdots,s_n,q,t)$ in Eq.(\ref{f}) is
an eigenfunction of the difference operator $D^1$
\begin{eqnarray}
D^1(s_1,s_2,\cdots,s_n,t,q) f(x_1,x_2,\cdots,x_n)=
\sum_{i=1}^n s_i \cdot f(x_1,x_2,\cdots,x_n). \label{difeq}
\end{eqnarray}
\end{con}

The case $n=2$ is easy, and will be treated in Section 5. 
When $n\geq 3$, however,
to prove Conjecture \ref{con} seems
a very complicated task, and it is an open problem.
We have checked it by a computer-aid calculation up to $n=5$ for small degrees in $x_i$'s.

In Section 6, 
a supporting argument for the case $n=3$ will be given.
We will prove that a consequence of Conjecture \ref{con}
(see  Eq. (\ref{raising}) below)
agrees with the theorem by Lassalle and Schlosser for the case $n=3$.

\section{Main Consequence}
We briefly recall the notion of the Macdonald polynomials \cite{Mac}.
Let $x_1,x_2,\cdots,x_n$ be a set of indeterminates, 
and $\Lambda_n={\bf Z}[x_1,\cdots,x_n]^{S_n}$ denotes the 
ring of symmetric polynomials. The ring of symmetric functions $\Lambda$
is defined as the inverse limit of the $\Lambda_n$
in the category of graded rings. 
Let $F={\bf Q}(q,t)$ be the field
of rational functions in independent indeterminates $q$ and $t$,
and set $\Lambda_{F}=\Lambda\otimes_{\bf Z} F$.

Let $p_n=\sum_i x_i^n$ be the power sum symmetric functions, and 
denote $p_\lambda=p_{\lambda_1}p_{\lambda_2}\cdots$
for any partition $\lambda=(\lambda_1,\lambda_2,\cdots)$.
The scalar product is introduced by
\begin{eqnarray}
\langle p_\lambda,p_\mu\rangle_{q,t}=
\delta_{\lambda,\mu}\prod_{i\geq1}i^{m_i}m_i! 
\prod_{j\geq 1} {1-q^{\lambda_j}\over 1-t^{\lambda_j}},\label{sp}
\end{eqnarray}
where $m_i=m_i(\lambda)$ is the 
multiplicity of the part $i$ in the partition $\lambda$.

The Macdonald symmetric functions $P_\lambda(x;q,t)\in \Lambda_F$ 
are uniquely characterized by
the following two conditions \cite{Mac}:
\begin{eqnarray}
\!\!\!\!\!\!\!\!\!\!\!\!
{\rm (a)}&&P_\lambda=m_\lambda+\sum_{\mu<\lambda} u_{\lambda\mu} m_\mu,
\end{eqnarray}
where $m_\lambda$ is the monomial symmetric function
associated with $\lambda$, $u_{\lambda\mu}\in F$, and the symbol ``$<$''
means the dominance ordering on the partitions.
\begin{eqnarray}
{\rm (b)}&&\langle P_\lambda,P_\mu\rangle_{q,t}=0\quad  {\rm if}\quad \lambda\neq \mu.
\end{eqnarray}

The dual of $P_\lambda$ with respect to the scalar product (\ref{sp})
is denoted by $Q_\lambda(x;q,t)=
b_\lambda(q,t)P_\lambda(x;q,t)$, where
$b_\lambda(q,t)=\langle P_\lambda,P_\lambda\rangle_{q,t}^{-1}$.
As for the explicit expression for $b_\lambda(q,t)$, see (6.19) of \cite{Mac}.

The symmetric function $g_n(x;q,t)\in \Lambda_F$ is defined by
\begin{eqnarray}
\prod_{i\geq 1}{(t x_i y;q)_\infty \over (x_i y;q)_\infty}=
\sum_{n\geq 0} g_n(x;q,t) y^n,
\end{eqnarray}
where $(a;q)_\infty =\prod_{i=0}^\infty (1-a q^i)$.
It is well known that we have $Q_{(n)}(x;q,t)=g_n(x;q,t)$
(equation (5.5) in \cite{Mac}). We use the convention that
$g_n(x;q,t)=0$ for $n<0$.
We write $g_a=g_{a_1}g_{a_2}\cdots g_{a_n}$
for any $a=(a_1,a_2,\cdots,a_n)\in {\bf Z}^n$.

Nextly, we recall the definition of the raising operators.
Let $a=(a_i,\cdots,a_n)\in {\bf Z}^n$. For each pair of integers $i,j$
satisfying $1\leq i<j\leq n$, define the action of $R_{ij}$ by
\begin{eqnarray}
R_{ij}(a)=(a_1,\cdots,a_i+1,\cdots,a_j-1,\cdots,a_n).
\end{eqnarray}
Any product of the form $\prod_{i<j}R_{ij}^{\theta_{i,j}}$ is called a raising operator.
For any raising operator $R$, $Rg_\lambda$ means $g_{R\lambda}$.
\bigskip

It was argued in the paper \cite{S2} that the solution to the equation
(\ref{difeq}) is interpreted as the raising operator for the
Macdonald polynomials. See Proposition A.6 in Appendix in \cite{S2}.
Therefore, the following is a consequence of Conjecture \ref{con}.
\begin{prop}\label{conse}
Let $\lambda=(\lambda_1,\cdots,\lambda_n)$ be a partition. 
Set $s_i=t^{n-i} q^{\lambda_i}$. 
Assume that Conjecture \ref{con} is true. Then 
the Macdonald polynomila $Q_\lambda=Q_\lambda(x;q,t)$
can be represented by the raising operator as
\begin{eqnarray}
Q_{(\lambda_1,\cdots,\lambda_n)}&=&
\prod_{1\leq k<\ell\leq n}
(1-R_{k\ell}) \label{raising}\\
&&\times
\sum_{(\theta)\in \mathsf{M}^{(n)}}
c_n(\{\theta_{i,j};1\leq i<j\leq n\};s_1,\cdots,s_n,q,t)
\prod_{1\leq i<j\leq n}
R_{ij}^{\theta_{i,j}}~
g_\lambda.\nonumber
\end{eqnarray}
\end{prop}
\bigskip

\section{Lassale and Schlosser's Theorem}
In this section, we recall 
Lassale and Schlosser's theorem for the raising operators of the
Macdonald polynomials \cite{LS} 
(which was announced in \cite{LS0}).

Let $u_1,\cdots, u_n$ be indeterminates, and $\theta_1,\cdots,\theta_n$
be nonnegative integres. Write $v_k=q^{\theta_k} u_k$ for simplicity.
Lassale and Schlosser obtained the following function
by inverting the Pieri formula
\begin{eqnarray}
&&C^{(q,t)}_{\theta_1,\cdots,\theta_n}(u_1,\cdots,u_n)
=\prod_{k=1}^n t^{\theta_k}
{(q/t;q)_{\theta_k} \over (q;q)_{\theta_k}}
{(qu_k;q)_{\theta_k} \over (qt u_k;q)_{\theta_k}}
\prod_{1\leq i<j\leq n}
{(qu_i/tu_j;q)_{\theta_i} \over (qu_i/u_j;q)_{\theta_i}}
{(tu_i/v_j;q)_{\theta_i} \over (u_i/v_j;q)_{\theta_i}}\nonumber\\
&&\qquad\qquad\qquad \times
{1\over \Delta(v)}
\mathop{{\rm det}}_{1\leq i,j\leq n}
\left[
v_i^{n-j}
\left(
1-
t^{j-1}
{1-tv_i\over 1-v_i}
\prod_{k=1}^n{u_k-v_i \over t u_k-v_i}
\right)
\right],\label{LSfun}
\end{eqnarray}
where $\Delta(v)=\prod_{1\leq i<j\leq n}(v_i-v_j)$.
For any $\theta =(\theta_{i,j})\in \mathsf{M}^{(n)}$, we write
\begin{eqnarray}
&&
\zeta_k(\theta)=\sum_{j=k+1}^n \theta_{k,j}-\sum_{j=1}^{k-1} \theta_{j,k},
\qquad 1\leq k\leq n,\\
&&
\xi_{ik}(\theta)=\sum_{j=k+2}^n (\theta_{i,j}-\theta_{k+1,j}),
\qquad 1\leq i\leq k\leq n-1.
\end{eqnarray}

The following important result was obtained
(Theorem 5.1 in \cite{LS}).
\begin{thm}[Lassale Schlosser]\label{LS}
Let $\lambda=(\lambda_1,\cdots,\lambda_n)$ be an arbitrary partition with length $n$.
We have
\begin{eqnarray}
\!\!\!\!&&Q_\lambda(q,t)\label{LStheorem}\\
\!\!\!\!&&=\sum_{\theta\in \mathsf{M}^{(n)}}
\prod_{k=1}^{n-1}
C^{(q,t)}_{\theta_{1,k+1},\cdots,\theta_{k,k+1}}
(\{u_i=q^{\lambda_i-\lambda_{k+1}+\xi_{ik}(\theta) }t^{k-i};1\leq i\leq k\})
\prod_{k=1}^n g_{\lambda_k+\zeta_k(\theta)}.\nonumber
\end{eqnarray}
\end{thm}

Comparing Lassalle and Schlosser's formula with our conjecture (\ref{raising}),
we observe the following.
\begin{con}\label{compare}
Setting $s_i=t^{n-i}q^{\lambda_i}$, we have
\begin{eqnarray}
\!\!\!\!\!&&\sum_{\theta\in \mathsf{M}^{(n)}}
\prod_{k=1}^{n-1}
C^{(q,t)}_{\theta_{1,k+1},\cdots,\theta_{k,k+1}}
(\{u_i=q^{\lambda_i-\lambda_{k+1}+\xi_{ik}(\theta) }t^{k-i};1\leq i\leq k\})
\prod_{1\leq i<j\leq n}R_{ij}^{\theta_{i,j}}\nonumber\\
\!\!\!\!\!&&=
\prod_{1\leq k<\ell \leq n}(1-R_{k\ell})\\
\!\!\!\!\!&&\times
\sum_{\theta\in \mathsf{M}^{(n)}}
c_n(\{\theta_{i,j};1\leq i<j\leq n\};
s_1,\cdots,s_n,q,t)\prod_{1\leq i<j\leq n}R_{ij}^{\theta_{i,j}}\nonumber,
\end{eqnarray}
where $R_{ij}$ denotes the raising operator.
\end{con}

One can prove Conjecture \ref{compare} for $n=2$ and $n=3$.
These two cases will be treated in the following two sections.
The case $n\geq 4$, however, is complicated and remains as an open problem.

It should be stressed that the identity in Conjecture \ref{compare} 
is a quite nontrivial one. 
In Section 6, an elementary proof for $n=3$ will be given.
It seems that some combinatorial properties for general $n$
hopefully can be extracted from that.
\medskip

In view of (\ref{c_n}), (\ref{LSfun}), (\ref{LStheorem}),  and $s_i=t^{n-i}q^{\lambda_i}$,
the difference between the two functions
\begin{eqnarray*}
&&\prod_{k=1}^{n-1}
C^{(q,t)}_{\theta_{1,k+1},\cdots,\theta_{k,k+1}}
(\{u_i=q^{\lambda_i-\lambda_{k+1}+\xi_{ik}(\theta) }t^{k-i};1\leq i\leq k\}),\\
&&c_n(\{\theta_{i,j};1\leq i<j\leq n\};
s_1,\cdots,s_n,q,t),
\end{eqnarray*}
only comes from the determinant factor in (\ref{LSfun}),
namely 
\begin{eqnarray}
&&{1\over \Delta(v)}
\mathop{{\rm det}}_{1\leq i,j\leq n}
\left[
v_i^{n-j}
\left(
1-
t^{j-1}
{1-tv_i\over 1-v_i}
\prod_{k=1}^n{u_k-v_i \over t u_k-v_i}
\right)
\right]\label{det}\\
&&=
\sum_{K\subset \{1,\cdots,n\}}(-1)^{|K|}(1/t)^{\left( |K|+1\atop 2\right)}
\prod_{k\in K\atop j\notin K}
{v_j-v_k/t\over v_j-v_k}
\prod_{i=1}^{n+1}\prod_{k\in K} {u_i-v_k \over u_i-v_k/t}.\nonumber
\end{eqnarray}
The RHS in the above equation is the formula (4.2) of the paper \cite{LS}.
In Section 6, we will use this expression for $n=1$ and $2$.

\section{Case $n=2$}
Let us consider the simplest case $n=2$. 
Setting $f(x_1,x_2;s_1,s_2,q,t)=(1-x_2/x_1)g(x_2/x_1)$,
$g(x)=\sum_{\theta=0}^\infty c_\theta x^\theta$, and $c_0=1$,
the difference equation (\ref{f}) for $n=2$ can be written as
\begin{eqnarray}
s_1(1-qt^{-1}x)g(q x)+s_2(1-q^{-1}t x)g(q^{-1}x)=
(s_1+s_2)(1-x)g(x).
\end{eqnarray}
Comparing the coefficients of $x^\theta$ from both sides,
the recurrence relation for the coefficients is obtained as
\begin{eqnarray}
c_\theta=t{(1-q^\theta t^{-1})(1-q^\theta t^{-1}s_1/s_2)\over 
(1-q^\theta)(1-q^\theta  s_1/s_2)}c_{\theta-1},
\qquad \theta=1,2,3,\cdots.
\end{eqnarray}
Hence we have 
\begin{eqnarray}
c_\theta=t^\theta{(qt^{-1};q)_\theta(qt^{-1}s_1/s_2;q)_\theta\over 
(q;q)_\theta(q s_1/s_2;q)_\theta},
\end{eqnarray}
which is given in  Eq. (\ref{c2}). Hence Conjecture \ref{con} is true for $n=2$.

Let us connect our series $f(x_1,x_2;s_1,s_2,q,t)$
with the formula of Jing and J\'ozefiak \cite{JJ}, and 
Lassale and Schlosser's one for $n=2$.
Note that the series can be rewritten as 
\begin{eqnarray}
f(x_1,x_2;s_1,s_2,q,t)
=(1-x_2/x_1)\sum_{\theta=0}^\infty
c_\theta
(x_2/x_1)^{\theta}\nonumber
=\sum_{\theta=0}^\infty
(c_\theta-c_{\theta-1})
(x_2/x_1)^{\theta},
\end{eqnarray}
if we set $c_{-1}=0$.
Then we observe that
\begin{eqnarray*}
&&c_\theta-c_{\theta-1}\\
&=&
t^{\theta} {(qt^{-1};q)_{\theta} \over (q;q)_{\theta} }
{(qt^{-1}s_1/s_2;q)_{\theta} \over (qs_1/s_2;q)_{\theta} }
\left(1-
t^{-1}{(1-q^\theta )(1-q^\theta s_1/s_2)\over 
(1-q^\theta t^{-1})(1-q^\theta t^{-1} s_1/s_2)}\right)\\
&=&
t^{\theta} {(qt^{-1};q)_{\theta} \over (q;q)_{\theta} }
{(qt^{-1}s_1/s_2;q)_{\theta} \over (qs_1/s_2;q)_{\theta} }
{(1-t^{-1} )(1-q^{2\theta} t^{-1}s_1/s_2)\over 
(1-q^\theta t^{-1})(1-q^\theta t^{-1} s_1/s_2)}\\
&=&
t^{\theta} {(t^{-1};q)_{\theta} \over (q;q)_{\theta} }
{(t^{-1}s_1/s_2;q)_{\theta} \over (qs_1/s_2;q)_{\theta} }
{1-q^{2\theta} t^{-1}s_1/s_2\over 
1-t^{-1} s_1/s_2}.
\end{eqnarray*}
Setting $s_1/s_2=t q^{\lambda_1-\lambda_2}$,
we recover Jing and J\'ozefiak's formula \cite{JJ} from the last line.
In view of (\ref{LSfun}) for $n=1$, one finds that
the second line corresponds to Lassalle and Schlosser's expression.
Namely, Conjecture \ref{compare} is true for $n=2$.

\section{Case $n=3$}
Next, let us examine the case $n=3$. Unfortunately, 
we do not have a method to solve the difference equation (\ref{f}) at this moment.
One may, however, prove that 
Conjecture \ref{compare} is true for $n=3$.
Therefore it is expected that Conjecture \ref{con} holds for $n=3$.

First, let us note the following identity
\begin{eqnarray}
&&(1-x_2/x_1)(1-x_3/x_1)(1-x_3/x_2)\nonumber\\
&&=1-{x_2\over x_1}-{x_3\over x_1}-{x_3\over x_2}+
{x_2\over x_1}{x_3\over x_1}+{x_2\over x_1}{x_3\over x_2}+{x_3\over x_1}{x_3\over x_2}-
{x_2\over x_1}{x_3\over x_1}{x_3\over x_2}\label{saseki}\\
&&=1-{x_2\over x_1}-{x_3\over x_2}
+\alpha\left({x_2\over x_1}{x_3\over x_2}-{x_3\over x_1}\right)+
\left({x_2\over x_1}\right)^2{x_3\over x_2}
+{x_3\over x_1}{x_3\over x_2}- {x_2\over x_1}{x_3\over x_1}{x_3\over x_2},\nonumber
\end{eqnarray}
with an arbitrary coefficient $\alpha$. From this we have 
\begin{lem}
The series $f(x_1,x_2,x_3;s_1,s_2,s_3,q,t)$ can be recast as 
\begin{eqnarray}
&&f(x_1,x_2,x_3;s_1,s_2,s_3,q,t)=
(1-x_2/x_1)(1-x_3/x_1)(1-x_3/x_2)\nonumber\\
&&\times\sum_{\theta \in \mathsf{M}^{(3)}}c(\theta_{1,2},\theta_{1,3},\theta_{2,3})
(x_2/x_1)^{\theta_{1,2}}(x_3/x_1)^{\theta_{1,3}}(x_3/x_2)^{\theta_{2,3}}\label{change}\\
&&=\sum_{\theta \in \mathsf{M}^{(3)}}
\widetilde{c}(\theta_{1,2},\theta_{1,3},\theta_{2,3})
(x_2/x_1)^{\theta_{1,2}}(x_3/x_1)^{\theta_{1,3}}(x_3/x_2)^{\theta_{2,3}},\nonumber
\end{eqnarray}
where we have denoted 
$c_3(\theta_{1,2},\theta_{1,3},\theta_{2,3};s_1,s_2,s_3,q,t)=
c(\theta_{1,2},\theta_{1,3},\theta_{2,3})$ for simplicity, and
\begin{eqnarray}
&&\widetilde{c}(\theta_{1,2},\theta_{1,3},\theta_{2,3})\nonumber\\
&&=c(\theta_{1,2},\theta_{1,3},\theta_{2,3})
-c(\theta_{1,2}-1,\theta_{1,3},\theta_{2,3})
-c(\theta_{1,2},\theta_{1,3},\theta_{2,3}-1)\nonumber\\
&&+\alpha(\theta_{1,2}-1,\theta_{1,3},\theta_{2,3}-1)c(\theta_{1,2}-1,\theta_{1,3},\theta_{2,3}-1)\\
&&
-\alpha(\theta_{1,2},\theta_{1,3}-1,\theta_{2,3})c(\theta_{1,2},\theta_{1,3}-1,\theta_{2,3})
\nonumber\\
&&+c(\theta_{1,2}-2,\theta_{1,3},\theta_{2,3}-1)+c(\theta_{1,2},\theta_{1,3}-1,\theta_{2,3}-1)
\nonumber\\
&&
-c(\theta_{1,2}-1,\theta_{1,3}-1,\theta_{2,3}-1),\nonumber
\end{eqnarray}
with arbitrary coefficients $\alpha(\theta_{1,2},\theta_{1,3},\theta_{2,3})$'s.
\end{lem}
\proof In view of (\ref{saseki}), we have
\begin{eqnarray*}
\!\!\!\!&&\mbox{RHS of (\ref{change})}\\
\!\!\!\!&&=
\sum_{\theta \in \mathsf{M}^{(3)}}
\left(
1-{x_2\over x_1}-{x_3\over x_2}
+\alpha(\theta_{1,2},\theta_{1,3},\theta_{2,3})
\left({x_2\over x_1}{x_3\over x_2}-{x_3\over x_1}\right)+
\left({x_2\over x_1}\right)^2{x_3\over x_2}
+{x_3\over x_1}{x_3\over x_2}- {x_2\over x_1}{x_3\over x_1}{x_3\over x_2}
\right)\\
\!\!\!\!&&\qquad \times
c(\theta_{1,2},\theta_{1,3},\theta_{2,3})
(x_2/x_1)^{\theta_{1,2}}(x_3/x_1)^{\theta_{1,3}}(x_3/x_2)^{\theta_{2,3}}\\
\!\!\!\!&&=\mbox{LHS of (\ref{change})},
\end{eqnarray*}
where we have assumed that $c(\theta_{1,2},\theta_{1,3},\theta_{2,3})=0$
if some of $\theta_{i,j}$'s are negative.
\qed

Our claim in this section is the following.
\begin{prop}
If we set 
\begin{eqnarray}
\alpha(\theta_{1,2},\theta_{1,3},\theta_{2,3})
&=&
{1-t^{-1}\over 1-q^{\theta_{1,2}}t^{-1}}
{1-q^{\theta_{1,3}-\theta_{2,3}}t s_1/s_2\over1-q^{\theta_{1,3}-\theta_{2,3}} s_1/s_2}
{1-q^{2\theta_{1,2}+\theta_{1,3}-\theta_{2,3}+1}t^{-1} s_1/s_2\over
1-q^{\theta_{1,2}+\theta_{1,3}-\theta_{2,3}+1} s_1/s_2},
\end{eqnarray}
the $\widetilde{c}(\theta_{1,2},\theta_{1,3},\theta_{2,3})$
is written as
\begin{eqnarray}
\widetilde{c}(\theta_{1,2},\theta_{1,3},\theta_{2,3}) =
C^{(q,t)}_{\theta_{1,2}}(q^{\theta_{1,3}-\theta_{2,3}} t^{-1}s_1/s_2)
C^{(q,t)}_{\theta_{1,3},\theta_{2,3}}(t^{-1} s_1/s_3,t^{-1} s_2/s_3).
\label{comp}
\end{eqnarray}
Here the RHS is
 Lassalle and Schlosser's function in Theorem \ref{LS} for $n=3$.
\end{prop}

\proof
Set 
\begin{eqnarray}
\beta(i,j,k)={c(\theta_{1,2}-i,\theta_{1,3}-j,\theta_{2,3}-k)
\over c(\theta_{1,2},\theta_{1,3},\theta_{2,3})},
\end{eqnarray}
for simplicity. Using (\ref{c3}) we have
\begin{eqnarray}
&&\beta(1,0,0)=
t^{-1}
{1-q^{\theta_{1,2}} \over 1-q^{\theta_{1,2}}t^{-1}}
{1-q^{\theta_{1,2}+\theta_{1,3}-\theta_{2,3}}s_1/s_2 \over
1-q^{\theta_{1,2}+\theta_{1,3}-\theta_{2,3}}t^{-1}s_1/s_2},\\
&&\beta(0,1,0)=
t^{-1}
{1-q^{\theta_{1,3}} \over 1-q^{\theta_{1,3}}t^{-1}}
{1-q^{\theta_{1,3}}s_1/s_3 \over
1-q^{\theta_{1,3}}t^{-1}s_1/s_3}\\
&&\times
{1-q^{\theta_{1,3}-\theta_{2,3}}t^{-1}s_1/s_2 \over 
1-q^{\theta_{1,3}-\theta_{2,3}}s_1/s_2}
{1-q^{\theta_{1,2}+\theta_{1,3}-\theta_{2,3}}s_1/s_2 \over 
1-q^{\theta_{1,2}+\theta_{1,3}-\theta_{2,3}}t^{-1}s_1/s_2}\nonumber\\
&&\times
{1-q^{\theta_{1,3}}s_1/s_2 \over 
1-q^{\theta_{1,3}}t^{-1}s_1/s_2}
{1-q^{\theta_{1,3}-\theta_{2,3}-1}s_1/s_2 \over 
1-q^{\theta_{1,3}-\theta_{2,3}-1}t s_1/s_2},\nonumber\\
&&\beta(0,0,1)=
t^{-1}
{1-q^{\theta_{2,3}} \over 1-q^{\theta_{2,3}}t^{-1}}
{1-q^{\theta_{2,3}}s_2/s_3 \over
1-q^{\theta_{2,3}}t^{-1}s_2/s_3}\\
&&\times
{1-q^{\theta_{1,3}-\theta_{2,3}+1}s_1/s_2 \over 
1-q^{\theta_{1,3}-\theta_{2,3}+1}t^{-1}s_1/s_2}
{1-q^{\theta_{1,2}+\theta_{1,3}-\theta_{2,3}+1}t^{-1}s_1/s_2 \over 
1-q^{\theta_{1,2}+\theta_{1,3}-\theta_{2,3}+1}s_1/s_2}\nonumber\\
&&\times
{1-q^{-\theta_{2,3}}s_1/s_2 \over 
1-q^{-\theta_{2,3}}t s_1/s_2}
{1-q^{\theta_{1,3}-\theta_{2,3}}t s_1/s_2 \over 
1-q^{\theta_{1,3}-\theta_{2,3}} s_1/s_2},\nonumber\\
&&\beta(1,0,1)=
t^{-2}
{1-q^{\theta_{1,2}} \over 1-q^{\theta_{1,2}}t^{-1}}
{1-q^{\theta_{1,3}-\theta_{2,3}+1}s_1/s_2 \over
1-q^{\theta_{1,3}-\theta_{2,3}+1}t^{-1}s_1/s_2}\\
&&\times
{1-q^{\theta_{2,3}} \over 1-q^{\theta_{2,3}}t^{-1}}
{1-q^{\theta_{2,3}}s_2/s_3 \over
1-q^{\theta_{2,3}}t^{-1}s_2/s_3}
{1-q^{-\theta_{2,3}}s_1/s_2 \over 
1-q^{-\theta_{2,3}}t s_1/s_2}
{1-q^{\theta_{1,3}-\theta_{2,3}}t s_1/s_2 \over 
1-q^{\theta_{1,3}-\theta_{2,3}} s_1/s_2},\nonumber\\
&&\beta(2,0,1)=
t^{-3}
{1-q^{\theta_{1,2}-1} \over 1-q^{\theta_{1,2}-1}t^{-1}}
{1-q^{\theta_{1,2}} \over 1-q^{\theta_{1,2}}t^{-1}}\\
&&\times
{1-q^{\theta_{1,3}-\theta_{2,3}+1}s_1/s_2 \over
1-q^{\theta_{1,3}-\theta_{2,3}+1}t^{-1}s_1/s_2}
{1-q^{\theta_{1,2}+\theta_{1,3}-\theta_{2,3}}s_1/s_2 \over 
1-q^{\theta_{1,2}+\theta_{1,3}-\theta_{2,3}}t^{-1}s_1/s_2}\nonumber\\
&&\times
{1-q^{\theta_{2,3}} \over 1-q^{\theta_{2,3}}t^{-1}}
{1-q^{\theta_{2,3}}s_2/s_3 \over
1-q^{\theta_{2,3}}t^{-1}s_2/s_3}
{1-q^{-\theta_{2,3}}s_1/s_2 \over 
1-q^{-\theta_{2,3}}t s_1/s_2}
{1-q^{\theta_{1,3}-\theta_{2,3}} t s_1/s_2 \over 
1-q^{\theta_{1,3}-\theta_{2,3}} s_1/s_2},\nonumber\\
&&\beta(0,1,1)=
t^{-2}
{1-q^{\theta_{1,3}} \over 1-q^{\theta_{1,3}}t^{-1}}
{1-q^{\theta_{1,3}}s_1/s_3 \over
1-q^{\theta_{1,3}}t^{-1}s_1/s_3}
{1-q^{\theta_{2,3}} \over 1-q^{\theta_{2,3}}t^{-1}}
{1-q^{\theta_{2,3}}s_2/s_3 \over
1-q^{\theta_{2,3}}t^{-1}s_2/s_3}\\
&&\times
{1-q^{\theta_{1,3}}s_1/s_2 \over 
1-q^{\theta_{1,3}}t^{-1}s_1/s_2}
{1-q^{-\theta_{2,3}}s_1/s_2 \over 
1-q^{-\theta_{2,3}}t s_1/s_2},\nonumber\\
&&\beta(1,1,1)=
t^{-3}
{1-q^{\theta_{1,2}} \over 1-q^{\theta_{1,2}}t^{-1}}
{1-q^{\theta_{1,2}+\theta_{1,3}-\theta_{2,3}}s_1/s_2 \over
1-q^{\theta_{1,2}+\theta_{1,3}-\theta_{2,3}}t^{-1}s_1/s_2}\\
&&\times
{1-q^{\theta_{1,3}} \over 1-q^{\theta_{1,3}}t^{-1}}
{1-q^{\theta_{1,3}}s_1/s_3 \over
1-q^{\theta_{1,3}}t^{-1}s_1/s_3}
{1-q^{\theta_{2,3}} \over 1-q^{\theta_{2,3}}t^{-1}}
{1-q^{\theta_{2,3}}s_2/s_3 \over
1-q^{\theta_{2,3}}t^{-1}s_2/s_3}\nonumber\\
&&\times
{1-q^{\theta_{1,3}}s_1/s_2 \over 
1-q^{\theta_{1,3}}t^{-1}s_1/s_2}
{1-q^{-\theta_{2,3}}s_1/s_2 \over 
1-q^{-\theta_{2,3}}t s_1/s_2}.\nonumber
\end{eqnarray}
Write
\begin{eqnarray}
&&a_{12}=
t^{-1}
{1-q^{\theta_{1,2}} \over 1-q^{\theta_{1,2}}t^{-1}}
{1-q^{\theta_{1,2}+\theta_{1,3}-\theta_{2,3}}s_1/s_2 \over
1-q^{\theta_{1,2}+\theta_{1,3}-\theta_{2,3}}t^{-1}s_1/s_2},\\
&&
a_{13}=
t^{-1}
{1-q^{\theta_{1,3}-\theta_{2,3}}t^{-1}s_1/s_2 \over 
1-q^{\theta_{1,3}-\theta_{2,3}}s_1/s_2}
{1-q^{\theta_{1,3}} \over 1-q^{\theta_{1,3}}t^{-1}}
{1-q^{\theta_{1,3}}s_1/s_2 \over 1-q^{\theta_{1,3}}t^{-1}s_1/s_2}
{1-q^{\theta_{1,3}}s_1/s_3 \over 1-q^{\theta_{1,3}}t^{-1}s_1/s_3},\\
&&
a_{23}=
t^{-1}
{1-q^{\theta_{1,3}-\theta_{2,3}}t s_1/s_2 \over 
1-q^{\theta_{1,3}-\theta_{2,3}}s_1/s_2}
{1-q^{\theta_{2,3}} \over 1-q^{\theta_{2,3}}t^{-1}}
{1-q^{-\theta_{2,3}}s_1/s_2 \over 1-q^{-\theta_{2,3}}t s_1/s_2}
{1-q^{\theta_{2,3}}s_2/s_3 \over 1-q^{\theta_{2,3}}t^{-1}s_2/s_3},\\
&&
a_{13,23}=
t^{-2}
{1-q^{\theta_{1,3}} \over 1-q^{\theta_{1,3}}t^{-1}}
{1-q^{\theta_{1,3}}s_1/s_2 \over 1-q^{\theta_{1,3}}t^{-1}s_1/s_2}
{1-q^{\theta_{1,3}}s_1/s_3 \over 1-q^{\theta_{1,3}}t^{-1}s_1/s_3}\nonumber\\
&&\qquad\qquad\times{1-q^{\theta_{2,3}} \over 1-q^{\theta_{2,3}}t^{-1}}
{1-q^{-\theta_{2,3}}s_1/s_2 \over 1-q^{-\theta_{2,3}}t s_1/s_2}
{1-q^{\theta_{2,3}}s_2/s_3 \over 1-q^{\theta_{2,3}}t^{-1}s_2/s_3},
\end{eqnarray}
for notational simplicity.
Then we have
\begin{eqnarray}
&&\beta(1,0,0)=a_{12},\qquad 
\beta(0,1,1)=a_{13,23},\qquad \beta(1,1,1)=a_{12}a_{13,23},\nonumber\\
&&\alpha(\theta_{1,2},\theta_{1,3}-1,\theta_{2,3})\beta(0,1,0)=
(1-a_{12})a_{13},\\
&&-\beta(0,0,1)+
\alpha(\theta_{1,2}-1,\theta_{1,3},\theta_{2,3}-1)\beta(1,0,1)+\beta(2,0,1)=
-(1-a_{12})a_{23}.\nonumber
\end{eqnarray}
Thus
\begin{eqnarray}
&&1
-\beta(1,0,0)
-\beta(0,0,1)\nonumber\\
&&+\alpha(\theta_{1,2}-1,\theta_{1,3},\theta_{2,3}-1)
\beta(1,0,1)
-\alpha(\theta_{1,2},\theta_{1,3}-1,\theta_{2,3})
\beta(0,1,0)
\nonumber\\
&&+\beta(2,0,1)+\beta(0,1,1)
-\beta(1,1,1)\nonumber\\
&&=
(1-a_{12})(1-a_{13}-a_{23}+a_{13,23}), \label{beta}
\end{eqnarray}
holds. By using (\ref{det}), one can check that
RHS of (\ref{beta}) is exactly the determinant factor from
 Lassalle and Schlosser's expression.  Namely we have
\begin{eqnarray}
\mbox{RHS of (\ref{beta})}=
{C^{(q,t)}_{\theta_{1,2}}(q^{\theta_{1,3}-\theta_{2,3}} t^{-1}s_1/s_2)
C^{(q,t)}_{\theta_{1,3},\theta_{2,3}}(t^{-1} s_1/s_3,t^{-1} s_2/s_3)\over 
c_3(\theta_{1,2},\theta_{1,3},\theta_{2,3};s_1,s_2,s_3,q,t) }.
\end{eqnarray}
This implies Eq. (\ref{comp}).
\qed

\section{Some Special Cases}
If $q=t$, the difference equation (\ref{f}) can be immediately solved for geneal $n$. 
Namely, we have
\begin{eqnarray}
&&D^1(s_1,\cdots,s_n,q,q)\prod_{1\leq i<j\leq n}(1-x_j/x_i)
=(s_1+\cdots+s_n)\prod_{1\leq i<j\leq n}(1-x_j/x_i).\nonumber
\end{eqnarray}
This means
that $c_n(\{\theta_{i,j};1\leq i<j\leq n\};s_1,\cdots,s_n,q,q)=0$ 
except if $\theta_{i,j}=0$ for all $i,j$.
Hence Conjecture \ref{con} is true for $q=t$.
Since we have $Q_\lambda(x;q,q)=s_\lambda(x)$ (Schur function), and 
$g_n(x;q,q)=h_n$ (complete symmetric function),
the Jacobi-Trudi formula for the Schur
polynomials (see formula (3.4) in \cite{Mac}) is recovered from 
our conjecture Eq. (\ref{raising})
\begin{eqnarray}
s_\lambda=\prod_{1\leq i<j\leq n}(1-R_{ij})h_{\lambda}=
\mathop{\rm det}_{1\leq i,j\leq n} ( h_{\lambda_i-i+j}).\label{Schur}
\end{eqnarray}

Next, let $k$ be a positive integer. For $t=q^k$, the coefficients 
$c_n(\{\theta_{i,j};1\leq i<j\leq n\};s_1,\cdots,s_n,q,q^k)$ vanish
if $\theta_{i,j}\geq k$ for some $i,j$, and the series (\ref{f})
becomes truncated. 
Therefore, the difference equation (\ref{difeq})
reduces to an identity of Laurent polynomials in $x_i$'s.
(Note the denominator in $D^1$ is cancelled by the factor $\prod_{j<i}
(1- x_i/x_j)$ in $f$.)
Even in this case, the equation (\ref{difeq}) is still complicated and
we are not able to prove (\ref{difeq}) at present.
We have proved, by a computer-aid calculation, that Conjecture \ref{con} is true for 
the cases: (1) $n=3$ and $t=q^2,q^3,q^4$, (2) $n=4$ and $t=q^2$.

Finally, we argue the case $q=0$. 
The $q=0$ limit of $c_n(\{\theta_{i,j};1\leq i<j\leq n\};s_1,\cdots,s_n,q,t)$
can be examined in several manners.
One may apply the automorphism $\omega_{q,t}$ 
(defined by  $\omega_{q,t}(p_r)=(-1)^{r-1}{1-q^r\over 1-t^r}p_r$) to (\ref{raising}),
and use the method presented in Section 7 of \cite{LS}. 
Even if we consider the limit $q=0$ in this way, 
it seems a difficult task to prove Conjecture \ref{compare}.
Instead of going in this direction, we give another argument
from which the $q=0$ limit can be studied. 

In \cite{S1,S2}, another type of conjecture
for the series satisfying (\ref{difeq}) was obtained for $n=3$.
Let us recall the statement.
\begin{con}\label{com3}
The series
\begin{eqnarray}
&&f(x_1,x_2,x_3;s_1,s_2,s_3,q,t)\nonumber\\
&=&
\sum_{k=0}^\infty
{
(q t^{-1};q)_k(q t^{-1};q)_k(t;q)_k(t;q)_k \over 
(q;q)_k(q s_1/s_2;q)_k(q s_2/s_3;q)_k
(q s_1/s_3;q)_k}
 (q s_1/s_3)^k (x_3/x_1)^k\label{g-fun}\\
&&\times
\prod_{1\leq i<j\leq 3}
(1-x_j/x_i)
\cdot {}_2\phi_1\left(
{  q^{k+1} t^{-1}, q t^{-1}s_i/s_j 
\atop
q^{{k+1}}s_i/s_j  };q, t x_j/x_i
\right),\nonumber
\end{eqnarray}
satisfies the difference equation (\ref{difeq}) for $n=3$.
\end{con}
Here we have used the standard notation for the 
basic hypergeometric series
\begin{eqnarray*}
{}_2\phi_1\left(
{  a, b 
\atop
c };q, x
\right)=\sum_{n=0}^\infty {(a;q)_n(b;q)_n\over (q;q)_n(c;q)_n}x^n.
\end{eqnarray*}

It has not been proved that the above series (\ref{g-fun}) 
and the one given by (\ref{f}) for $n=3$ are the same.
One can check the agreement up to certain degree in $x_i$'s,
and we observe the following
\begin{con}The identity
\begin{eqnarray}
&&\sum_{\theta\in \mathsf{M}^{(3)}}
t^{\theta_{1,2}} {(qt^{-1};q)_{\theta_{1,2}} \over (q;q)_{\theta_{1,2}} }
{(q^{\theta_{1,3}-\theta_{2,3}}qt^{-1}s_1/s_2;q)_{\theta_{1,2}} \over 
(q^{\theta_{1,3}-\theta_{2,3}}qs_1/s_2;q)_{\theta_{1,2}} }\\
&&\times
t^{\theta_{1,3}} {(qt^{-1};q)_{\theta_{1,3}} \over (q;q)_{\theta_{1,3}} }
{(qt^{-1}s_1/s_3;q)_{\theta_{1,3}} \over 
(qs_1/s_3;q)_{\theta_{1,3}} }
t^{\theta_{2,3}} {(qt^{-1};q)_{\theta_{2,3}} \over (q;q)_{\theta_{2,3}} }
{(qt^{-1}s_2/s_3;q)_{\theta_{2,3}} \over 
(qs_2/s_3;q)_{\theta_{2,3}} }\nonumber\\
&&\times
{(qt^{-1}s_1/s_2;q)_{\theta_{1,3}} \over 
(qs_1/s_2;q)_{\theta_{1,3}} }
{(q^{-\theta_{2,3}}ts_1/s_2;q)_{\theta_{1,3}} \over 
(q^{-\theta_{2,3}}s_1/s_2;q)_{\theta_{1,3}} }
(x_2/x_1)^{\theta_{1,2}}(x_3/x_1)^{\theta_{1,3}}(x_3/x_2)^{\theta_{2,3}}\nonumber\\
&&
=\sum_{k=0}^\infty
{
(q t^{-1};q)_k(q t^{-1};q)_k(t;q)_k(t;q)_k \over 
(q;q)_k(q s_1/s_2;q)_k(q s_2/s_3;q)_k
(q s_1/s_3;q)_k}
 (q s_1/s_3)^k (x_3/x_1)^k\nonumber\\
&&\times
\prod_{1\leq i<j\leq 3} {}_2\phi_1\left(
{  q^{k+1} t^{-1}, q t^{-1}s_i/s_j 
\atop
q^{{k+1}}s_i/s_j  };q, t x_j/x_i
\right),\nonumber
\end{eqnarray}
holds.
\end{con}

Now an explanation about the author's heuristic argument is in order.  Note that 
the factor $\prod_{1\leq i<j\leq 3}(1-x_j/x_i)$ can be seen 
in the series (\ref{g-fun}). From this, one may expect that
the same factor can be factored out from Lassale and Schlosser's
expression (\ref{LStheorem}).
Assuming this factorization, one can arrive at the series (\ref{f}) after some exploration.
\medskip

Note that for the case $q=t$, the raising operator 
formula for the Schur polynomials (\ref{Schur})
is correctly derived from (\ref{g-fun}).  Since the RHS vanishes
except for $k=0$, and 
\begin{eqnarray*}
{}_2\phi_1\left(
{  1, s_i/s_j
\atop
q s_i/s_j };q, q x_j/x_i
\right)=1,
\end{eqnarray*} 
we have
\begin{eqnarray}
&&f(x_1,x_2,x_3;s_1,s_2,s_3,q,q)=
\prod_{1\leq i<j\leq 3}
{(1-x_j/x_i)},\nonumber
\end{eqnarray}
  from (\ref{g-fun}).

Assume Conjecture \ref{com3}, and 
consider the $q=0$ limit from (\ref{g-fun}).
Since we have set $s_i=t^{n-i} q^{\lambda_i}$ for the partition 
$\lambda=(\lambda_1,\cdots,\lambda_n)$, 
$\lim_{q\rightarrow 0} q s_i/s_j =0$ holds for $i<j$.
It can be seen that the RHS of (\ref{g-fun}) vanishes except for $k=0$, and
\begin{eqnarray}
 \lim_{q\rightarrow 0}{}_2\phi_1\left(
{  q t^{-1}, q t^{-1}s_i/s_j 
\atop
q s_i/s_j  };q, t x_j/x_i
\right)=\sum_{n=0}^\infty (t x_j/x_i)^n={1\over 1-t x_j/x_i}.
\end{eqnarray}
Namely we have
\begin{eqnarray}
&&f(x_1,x_2,x_3;s_1,s_2,s_3,0,t)=
\prod_{1\leq i<j\leq 3}
{1-x_j/x_i\over 1-t x_j/x_i},\nonumber
\end{eqnarray}
from (\ref{g-fun}).
It is well known that $Q_\lambda(x;0,t)=Q_\lambda(x;t)$
(Hall-Littlewood symmetric function), and 
$g_n(x;0,t)=q_n(x;t)$, where 
\begin{eqnarray}
\prod_{i\geq 1}{1-t x_i y \over 1-x_i y}=
\sum_{n\geq 0} q_n(x;t) y^n.
\end{eqnarray}
Thus the raising operator expression for the
Hall-Littlewood functions with partitions with length three
\begin{eqnarray}
Q_\lambda(t)=\prod_{1\leq i<j\leq 3}
{1-R_{ij}\over 1-t R_{ij}} q_\lambda.
\end{eqnarray}
is recovered from (\ref{g-fun}) (see equation ($2.15'$) in \cite{Mac}).
\bigskip

\noindent
{\it Acknowledgment.}~~~The author is indebted to the referee, whose 
suggestions have improved the presentation of the paper. He thanks 
M. Lassalle for stimulating discussion and valuable comments. 
He is grateful to M. Noumi for having discussion and kind hospitality 
at Kobe Univ.
This work is supported by the  Grant-in-Aid for Scientific Research 
(C) 16540183.

\end{document}